\documentstyle[A4,12pt,amssymb,amsfonts]{article}
\begin{document}
\author{Nicusor Dan\footnote{Travail r\'ealis\'e avec le support du
contrat Cex05-D11-11/2005.}}
\title{Sur la conjecture de Zagier pour $n=4$}
\date{}
\maketitle
\parindent=0pt
\parskip=5pt
\newtheorem{theor}{Th\'eor\`eme}
\newtheorem{prop}[theor]{Proposition}
\newtheorem{cor}[theor]{Corollaire}
\newtheorem{lemma}[theor]{Lemme}
\newtheorem{defin}[theor]{D\'efinition}
\newtheorem{conj}{Conjecture}
\newtheorem{rem}[theor]{Remarque}
\newtheorem{nota}[theor]{Notations}
\def\vs{\vskip 4mm}
\def\gl{G^{\lambda}}
\def\ol{\Omega^\lambda}
\def\dl{{\Delta}^{\lambda}}
\def\gly{G^{\lambda}_Y}
\def\oly{\Omega^{\lambda}_Y}
\def\dly{{\Delta}^{\lambda}_Y}
\def\l{\lambda}
\def\eq{equation}
\def\lu{{\lambda}_1}
\def\ld{{\lambda}_2}
\def\lt{{\lambda}_3}
\def\glu{G^{\lu}}
\def\gld{G^{\ld}}
\def\glt{G^{\lt}}
\def\olu{\Omega^\lu}
\def\old{\Omega^\ld}
\def\olt{\Omega^\lt}
\def\dlu{{\Delta}^{\lu}}
\def\dld{{\Delta}^{\ld}}
 \def\dlt{{\Delta}^{\lt}}
\def\ddl{\frac{\partial}{\partial\l}}
\def\ddlu{\frac{\partial}{\partial\lu}}
\def\ddld{\frac{\partial}{\partial\ld}}
\def\ra{\rightarrow}
\def\ci{{\ C}^\infty}
\def\tix{\tilde{X}}
\def\pii{\pi^{-1}(y)}
\def\curg{courant de Green }
\def\fd{forme diff{\'e}rentielle }
\def\pix{\pi:\tix \ra X}
\def\sl{||s||^{2\l}}
\def\piy{\pi^{-1}(Y)}
\def\espc{espace analytique complexe de dimension finie }
\def\va{\varphi}
\def\ep{\epsilon}
\def\pa{\partial}
\def\al{\alpha}
\def\n{\nabla}
\def\m{\nabla'}
\def\p{\nabla''}
\def\q{\nabla^2}
\def\pc{\prod_k}
\def\mn{\{1,\cdots ,m+n\}}
\def\lz{(\l_l)_l =0}
\def\bw{\bigwedge}
\def\be{\beta}
\def\cd{\cdot}
\def\cds{\cdots}
\def\miu{\mu}
\def\ga{\gamma}
\def\ro{\rho}
\def\xt{X\times T}
\def\de{\delta}
\def\om{\omega}
\def\ul{U_{(\l_l)}}
\def\ri{\rightline{$\Box$}}
\def\rid{\rightline{$\Box$}}
\def\dx{{\cal D}(X)}
\def\dxx{{\cal D}'(X)}
\def\Om{\Omega}
\def\ho{{\cal H}(\Om)}
\def\dih{{\cal D}'_h(X,\Om)}
\def\R{{\Bbb{R} } }
\def\Z{{\Bbb{Z} } }
\def\Q{{\Bbb{Q} } }
\def\P{{\Bbb{P} } }
\def\N{{\Bbb{N} } }
\def\O{{\cal O}}
\def\D{{\cal D}}
\def\A{{\cal A}}
\def\E{{\cal E}}
\def\C{{\Bbb{C}}}
\def\widebar{\overline}
\def\exq{{E_{\Q}^{\times}}}
\def\fxq{{F_{\Q}^{\times}}}

\begin{abstract}
We express a general $4$-hyperlogarithm as a linear combination of
$4$-hyperlogarithms in two variables. We reduce the Zagier's
conjecture for $n=4$ to a combinatorial statement. We give a short
survey of the strategy of Goncharov and Zagier for reducing the
Zagier's conjecture for general $n$ to combinatorial relations
between hyperlogarithms. Such a survey is missing in the
literature.

\vs

\ \ \ \ \ \ \ \ \ \ \ \ \ \ \ \ \ \ \ \ \ \ \ \ \ \ \ \ \ \ \ \ \
\ \ \ \ \ {\bf R\'esum\'e} \vs

On exprime un $4$-hyperlogarithme g\'en\'eral comme combinaison
lin\'eaire de $4$-hyperlogarithmes en deux variables. On r\'eduit
la conjecture de Zagier pour $n=4$ \`a un \'enonc\'e combinatoire.
On donne une pr\'esentation synth\'etique de la strat\'egie de
Goncharov et Zagier pour la r\'eduction de la conjecture de Zagier
pour $n$ g\'en\'eral \`a des relations combinatoires entre
hyperlogarithmes. Une telle synth\`ese n'existe pas dans la
litt\'erature.

\end{abstract}

\section{\'Enonc\'e de la conjecture}
\label{100}

\vs Soit $n$ un entier positif. Le polylogarithme de poids $n$, ou
le $n$-logarithme, est la fonction complexe d\'efinie sur le
disque unit\'e $|z|\leq 1$ par la s\'erie absolument convergente
$$P_n(z) = Li_n(z)= \sum_{k=1}^{\infty}\frac{z^k}{k^n}.$$
On observe que $Li_1(z) = - {\mathrm log} (1-z) = \int_0^z
\frac{dt}{1-t}$ et que

\begin{\eq}
\label{0500}  Li_n(z)= \int_0^z Li_{n-1}(z)\frac{dt}{t}
\end{\eq}

On d\'eduit que la fonction $P_n(z)$ se prolonge analytiquement
\`a une fonction holomorphe multivalu\'e sur $\C \setminus \{ 0,1
\} $. On lui attache ([Z1]) la fonction r\'eelle univalu\'ee,
continue sur $\C$ et analytiquement r\'eelle sur $\C \setminus \{
0,1 \}  $
$$\R P_n (z) = {\cal R}_n (\sum_{k=0}^{m-1} \frac{B_k}{k!}
{\mathrm log}^k(z\bar{z}) Li_{n-k}(z)),$$

o\`u $B_k$ sont les nombres de Bernoulli et o\`u ${\cal R}_n$
d\'esigne la partie r\'eelle si $n$ est impair et la partie
imaginaire si $n$ est pair. On trouve dans [BD] une
interpr\'etation en th\'eorie de Hodge de la fonction $\R P_n$.

\vs Pour tout corps $E$, on note $\Q [E\setminus\{ 0,1 \}  ]$
l'espace vectoriel sur $\Q$ ayant comme base les symboles $[x]$
pour chaque $x \in E\setminus\{ 0,1 \}$. On prolonge par
lin\'earit\'e la fonction $\R P_n$ \`a une application lin\'eaire
$\R P_n :\Q [\C\setminus\{ 0,1 \}] \to \R$. La conjecture de
Zagier ([Z1]) est la suivante:
\begin{conj} Soit $F$ un corps de nombres. Soient $\sigma_{r_2
+1},\cdots,\sigma_{r_1 + r_2}:F \to \R$ ses plongements r\'eels et
$\sigma_1=\widebar{\sigma_{r_1 + r_2 +
1}},\cdots,\sigma_{r_2}=\widebar{\sigma_{r_1 + 2r_2}}:F \to \C$
ses plongements complexes. Soit $n\geq 2$ un entier. On note $d_n
= r_1 + r_2$ si $n$ est impair et $d_n = r_2$ si $n$ est pair.
Alors il existe des \'el\'ements $y_1, \cdots, y_{d_n} \in \Q
[F\setminus\{ 0,1 \} ]$ tels que
$$\zeta_F(n)= \pi^{(r_1+2r_2 -d_n)n} |D_F|^{-1/2} {\mathrm det}(\R P_n
(\sigma_i(y_j))) \ \ \ \ \ \ \ \ 1\leq i,j \leq d_n,$$ o\`u
$\zeta_F(s)$ est la fonction z\^eta de Dedekind et $D_F$ le
discriminant de $F$. \label{conj1}
\end{conj}
En fait, comme on verra dans sa formulation en K-th\'eorie
alg\'ebrique (Conjecture $\ref{conj2}$, point c)), la conjecture a
une forme plus pr\'ecise, dans laquelle les \'el\'ements $y_1,
\cdots, y_{d_n}$ sont cycles dans un certain sens.

Pour $F=\Q$ et $n$ impair la conjecture est triviale, car
$\zeta_{\Q}(n)=\sum_{k=1}^{\infty} \frac{1}{k^n}=\R P_n (1)$. Pour
$F$ totalement r\'eel et $n$ pair, la conjecture \'equivaut \`a
$\zeta_F(n) \in \pi^{r_1 n} \Q$. Pour $F$ g\'en\'eral et $n=1$ le
polylogarithme est le logarithme classique et si on consid\`ere
dans l'\'enonc\'e le r\'esidu $\zeta_F^{\times}(1)$ au lieu de
$\zeta_F(1)$, la conjecture est une variante faible du
th\'eor\`eme de Dedekind. Pour $F$ g\'en\'eral, si $n=2$ la
conjecture est un th\'eor\`eme de Zagier ([Z2]) et si $n=3$ un
th\'eor\`eme de Goncharov ([G1]).

\section{La strat\'egie pour prouver la conjecture de Zagier (d'apr\`es Goncharov et Zagier)}
\label{200} Le th\'eor\`eme $\ref{theor2}$ est le seul r\'esultat
originel de cette section.
\subsection{Pas 1: R\'eduction \`a un \'enonc\'e de K-th\'eorie
alg\'ebrique} \label{2.1.} Pour chaque corps $E$, suivant Zagier,
on d\'efinit par induction sur $n$ le sous-espace vectoriel ${\cal
R}_n^{\cal P}(E)\subset \Q [E\setminus\{ 0,1 \}  ]$ des "relations
entre polylogarithmes on $E$" et on pose ${\cal P}_n(E):= \Q
[E\setminus\{ 0,1 \}]/{\cal R}_n^{\cal P}(E)$. On note $[x]_n$ la
classe de $[x]$ modulo ${\cal R}_n^{\cal P}(E)$. Le sous-espace
vectoriel ${\cal R}_1^{\cal P}(E)$ est par d\'efinition engendr\'e
par $[x]-[y]-[x/y]$, pour $x,y\in E\setminus\{ 0,1 \} , x\neq y.$
Donc ${\cal P}_1(E)=E_{\Q}^{\times}$. On consid\`ere le morphisme
$\delta_2 : \Q [E\setminus\{ 0,1 \}] \to \Lambda^2
E_{\Q}^{\times}$ donn\'e par $[x]\to (1-x)\wedge x$ et les
morphismes $\delta_n :\Q [E\setminus\{ 0,1 \}] \to {\cal
P}_{n-1}(E) \otimes E_{\Q}^{\times}$ donn\'es par $[x]\to
[x]_{n-1} \otimes x$ si $n\geq 3$. On note ${\cal K}_n(E)$ le
noyau ${\mathrm Ker} \delta_n$. On d\'efinit ${\cal R}_n^{\cal
P}(E)$ comme le sous-espace vectoriel engendr\'e par $\alpha(1) -
\alpha(0)$ pour tous les \'el\'ements $\alpha$ de ${\cal
K}_n(E(t)),$ o\`u $t$ est une variable. On prouve que $\delta_n
({\cal R}_n^{\cal P}(E))=0$ et on obtient donc des applications
$$\delta_2 : {\cal P}_2(E)\to \Lambda^2
E_{\Q}^{\times}$$
\begin{\eq}
\label{1000} \ \ \ \ \ \ \ \ \ \ \ \ \ \ \ \ \ \ \ \ \ \ \ \ \ \ \
\ \ \ \ \delta_n : {\cal P}_n(E)\to {\cal P}_{n-1}(E) \otimes
E_{\Q}^{\times}\ \ \ \ \ \ \ \ \ \ \ (n\geq 3)
\end{\eq}
On prouve que $\R P_n ({\cal R}_n^{\cal P}(\C))=0$ et on obtient
donc une application $\R P_n : {\cal P}_n(\C)\to \R$. La
formulation en K-th\'eorie alg\'ebrique de la conjecture de Zagier
est la suivante:
\begin{conj}
\label{conj2} Pour chaque corps $E$, il existe une application
$r_{{\cal P}_n} : K_{2n-1}(E)_{\Q} \to {\cal P}_n(E)$ qui: a)
v\'erifie $\R P_n \circ  r_{{\cal P}_n} =r_n$ si $E=\C$; b) est
naturelle pour les inclusions de corps; c) a l'image dans
${\mathrm Ker} \delta_n \subset {\cal P}_n(E)$.
\end{conj}

Cette conjecture est un cas particulier d'une conjecture plus
optimiste, qui dit que le complexe
\begin{\eq}
\label{2000} {\cal P}_n(E) \stackrel{\delta_n}{\longrightarrow}
{\cal P}_{n-1}(E) \otimes \exq
\stackrel{\delta_{n-1}}{\longrightarrow} {\cal P}_{n-2}(E) \otimes
\Lambda^2 \exq  \stackrel{\delta_{n-2}}{\longrightarrow}\cdots
\stackrel{\delta_{2}}{\longrightarrow}\Lambda^n \exq
\end{\eq}
est le complexe motivique $\underline{\Q(n)}$ sur $ {\mathrm Spec}
E$ (donc ${\mathrm Ker} \delta_n = K_{2n-1}^{[n]}(E)_{\Q}$, le
facteur de $K_{2n-1}(E)_{\Q}$ de poids $n$ pour les op\'erations
d'Adams).

\begin{theor}
\label{theor1}
La conjecture $\ref{conj2}$ implique la conjecture
$\ref{conj1}$.
\end{theor}

\par{\bf Preuve: }
Soit $B_{{\mathrm ct}}GL(\C)$ le classifiant du groupe topologique
$GL(\C)$ et $B_{{\mathrm dis}}GL(\C)$ le classifiant du groupe
$GL(\C)$ vu comme groupe discret. Il existe un \'el\'ement
canonique $b_{2n-1}\in H^{2n-1} (B_{{\mathrm ct}}GL(\C),\R)$. Il
induit un \'el\'ement dans $H^{2n-1} (B_{{\mathrm
dis}}GL(\C),\R)$, donc un morphisme $H_{2n-1}(B_{{\mathrm
dis}}GL(\C),\R)\to \R$. On compose ce morphisme avec le morphisme
de Hurewicz
$$K_{2n-1}:=\pi_{2n-1}(B_{{\mathrm dis}}GL(\C)^{+})\to
H_{2n-1}(B_{{\mathrm dis}}GL(\C)^{+},\Z)=H_{2n-1}(B_{{\mathrm
dis}}GL(\C),\Z)$$ et on obtient un morphisme $r_n:K_{2n-1}(\C) \to
\R$. Pour tout groupe ab\'elien $A$, on note
$A_{\Q}:=A\otimes_{\Z}{\Q}$. Le th\'eor\`eme de Borel ([B]) dit
que $K_{2n-1}(F)$ est un groupe ab\'elien de rang $d_n$ et que, si
on fixe une base $x_1,\cdots,x_{d_n}$ de $K_{2n-1}(F)_{\Q}$, on a
$$\zeta_F(n)\in {\Q}^{\times} \pi^{(r_1+2r_2 -d_n)n} |D_F|^{-1/2} {\mathrm det}(r_n
(\sigma_i(y_j))) \ \ \ \ \ \ \ \ 1\leq i,j \leq d_n.$$ Il suffit
de prendre $y_i = r_{{\cal P}_n} (x_i)$ pour $i=1,\cdots,d_n$.

\subsection{Pas 2: R\'eduction du r\'egulateur \`a un hyperlogarithme
(ou polylogarithme grassmannien, ou polylogarithme d'Aomoto)}
\label{2.2.}

Si $a$ et $b$ sont deux points sur une vari\'et\'e complexe $X$ et
$\omega_1, \cdots, \omega_n$ sont des $1-$formes holomorphes sur
$X$, l'int\'egrale it\'er\'ee se d\'efinit par induction sur $n$
par la formule
$$\int_a^b \omega_1 \circ  \cdots \circ  \omega_n = \int_a^b (\int_a^t \omega_1 \circ
\cdots\circ  \omega_{n-1})\omega_n(t).$$

On peut \'ecrire la formule $(\ref{0500})$ comme
$$Li_n(z)=\int_0^z \frac{dt}{1-t} \circ  \frac{dt}{t} \circ  \cdots \circ
\frac{dt}{t}.$$

En g\'en\'eralisant cette formule, on d\'efinit les
hyperlogarithmes comme les int\'egrales it\'er\'ees
$$H(a_0|a_1,\cdots,a_n|a_{n+1})=\int_{a_0}^{a_{n+1}} \frac{dt}{t-a_1} \circ  \frac{dt}{t-a_2} \circ  \cdots \circ
\frac{dt}{t-a_n}.$$

C'est une fonction complexe multivalu\'ee sur l'ensemble des
$(n+2)-$ulpes complexes $(a_0,\cdots,a_{n+1})$ v\'erifiant $a_0
\neq a_1, a_n \neq a_{n+1}$ (pour que l'int\'egrale converge). On
peut lui associ\'er une fonction univalu\'ee r\'eelle $\R
H(a_0|a_1,\cdots,a_n|a_{n+1})$.

On peut g\'en\'eralis\'er la construction du paragraphe
$\ref{2.1.}$. On note $E_{\times}^{n+2}$ l'ensemble des
$(n+2)-$uples $(a_0,\cdots,a_{n+1})$ de $E$ satisfaisant $a_0\neq
a_1$ et $a_n\neq a_{n+1}$. On note $\Q [ E_{\times}^{n+2}]$
l'\'espace vectoriel sur $\Q$ ayant comme base les symboles
$[a_0|a_1, \cdots, a_n|a_{n+1}]$ pour $(a_0, \cdots, a_{n+1})\in
E_{\times}^{n+2}$. On d\'efinit par induction l'espace vectoriel
${\cal R}_n^{\cal H}(E)\subset \Q [ E_{\times}^{n+2}]$ des
"relations entre hyperlogarithmes on $E$" et on pose ${\cal
H}_n(E):=\Q [ E_{\times}^{n+2}]/{\cal R}_n^{\cal H}(E)$. On note
toujours $[a_0|a_1, \cdots, a_n|a_{n+1}]$ la classe de
l'\'el\'ement $[a_0|a_1, \cdots, a_n|a_{n+1}]$ modulo ${\cal
R}_n^{\cal H}(E)$. Les morphismes canoniques ${\cal P}_n(E) \to
{\cal H}_n(E)$ sont des isomorphismes pour $n\leq 3$ et seulement
des injections pour $n>3$. On a des morphismes
\begin{\eq}
\label{3000} \delta_n : {\cal H}_n(E) \to \oplus_{0<k<n/2} [{\cal
H}_{n-k}(E) \otimes {\cal H}_k(E)] \oplus \Lambda^2{\cal
H}_{n/2}(E)
\end{\eq}
(le dernier terme dispara\^it si $n$ est impair) et une
application "r\'ealisation"
\begin{\eq}
\label{4000}
\R {\cal H}_n :{\cal H}_n (\C) \to \R
\end{\eq}
Il existent deux autres g\'en\'eralisation des polylogarithmes,
qu'on ne d\'efinit pas ici. Les polylogarithmes grassmanniens sont
param\'etr\'es par les sous-espaces vectoriels de dimension $n$ de
$E^{2n}$, transverses aux hyperplans de coordon\'ees. Ils
d\'ependent de  $n^2$ variables. On leur attache des espaces
vectoriels ${\cal G}_n(E)$ comme ci-dessus v\'erifiant ${\cal
G}_1(E)=\exq$ et des applications $\delta_n$ comme $(\ref{3000})$
et $\R G_n:{\cal G}_n(\C)\to \R$ comme $(\ref{4000})$. Les
polylogarithmes d'Aomoto sont param\'etr\'es par $2n+2$ hyperplans
$(L_0,\cdots,L_n;M_0, \cdots,M_n)$ en ${\P}^n(E)$, satisfaisant
une condition de transversalit\'e, modulo l'action du groupe
$PGL(n+1)$. Ils d\'ependent de $n^2$ variables. On leur attache
des espaces vectoriels ${\cal A}_n(E)$ et des applications
$\delta_n, \R A_n$ comme ci-dessus. Toutes les conjectures qu'on
va formuler pour ${\cal H}_n(E)$ ont une variante identique avec
${\cal H}$ remplac\'e par ${\cal G}$ ou ${\cal A}$. On conjecture
en fait ${\cal H}_n(E)={\cal G}_n(E)={\cal A}_n(E)$.
\begin{theor}
\label{theor2} Pour chaque corps $E$, il existe une application
$r_{{\cal H}_n} : K_{2n-1}(E)_{\Q} \to {\cal H}_n(E)$ qui: a)
v\'erifie $\R H_n \circ  r_{{\cal P}_n} =r_n$ si $E=\C$; b) est
naturelle pour les inclusions de corps.
\end{theor}
\par{\bf Preuve: }
Le r\'egulateur de Borel co\"incide avec le r\'egulateur de
Beilinson multipli\'e par un rationnel non-nul. Goncharov
construit explicitement dans [G2] et [G3] le r\'egulateur de
Beilinson, en le factorisant par le complexe des polylogarithmes
grassmanniens. Donc le th\'eor\`eme est vrai si on remplace ${\cal
H}_n(E)$ par ${\cal G}_n(E)$.

Dans le paragraphe $7.3.$ de [G3], Goncharov donne la
r\'ealisation motivique des polylogarithmes grassmanniens. On
observe qu'elle s'assemble dans une famille de r\'ealisations
param\'etr\'es par un ouvert du grassmannien $G_{2n}^n$ des
sous-espaces vectoriels de $E^{2n}$ de dimension $n$. Le
th\'eor\`eme $5.6.$ de [G4] s'applique dans ce cas et implique que
le polylogarithme grassmannien est une combinaison lin\'eaire des
hyperlogarithmes. On d\'eduit une application $r_{{\cal G}_n {\cal
H}_n}: {\cal G}_n(E)\to {\cal H}_n(E)$ naturelle pour les
inclusions de corps et v\'erifiant $\R H_n \circ  r_{{\cal G}_n
{\cal H}_n}=\R G_n$ si $E=\C$, d'o\`u le th\'eor\`eme.

Il serait int\'eressant de d\'ecrire explicitement le morphisme
$r_{{\cal G}_n {\cal H}_n}$. A notre connaissance, personne n'a
fait cet exercice.

\begin{conj}
\label{conj3} Il existe une application $r_{{\cal H}_n}$ comme
dans le th\'eor\`eme $\ref{theor2}$ qui, en plus, a l'image
contenue dans ${\mathrm Ker} \delta_n$.
\end{conj}
La description explicite du morphisme $r_{{\cal G}_n {\cal H}_n}$
pourra aider a prouver cette conjecture.

\subsection{R\'eduction du hyperlogarithme (ou polylogarithme grassmannien,
ou polylogarithme d'Aomoto) \`a un polylogarithme}
\begin{conj}
\label{conj4} Il existe une application $r_{{\cal H}_n {\cal
P}_n}: ({\mathrm Ker} \delta_n \subset {\cal H}_n(E)) \to
({\mathrm Ker} \delta_n \subset {\cal P}_n(E))$ qui: a) v\'erifie
$\R P_n \circ r_{{\cal H}_n {\cal P}_n}=\R H_n$ si $E=\C$; b) est
naturelle pour les inclusions de corps.
\end{conj}
Si on \'ecrit ${\cal H}(E)=\oplus_{n=1}^{\infty} {\cal H}_n(E)$,
gradu\'e par $n\in \N$, les applications $(\ref{3000})$
d\'efinissent une cod\'erivation gradu\'ee $\delta:{\cal H}(E)\to
\Lambda^2 {\cal H}(E)$. La conjecture $\ref{conj4}$ est un cas
particulier d'une conjecture plus optimiste, qui affirme que le
sous-complexe de graduation $n$ du complexe
$$0\to {\cal H}\stackrel{\delta}{\longrightarrow}\Lambda^2{\cal H}\stackrel{\delta}{\longrightarrow}\cdots \stackrel{\delta}{\longrightarrow} \Lambda^n{\cal H}
\to 0$$ est quasi-isomorphe au complexe $(\ref{2000})$, d'une
mani\`ere compatible avec les inclusions de corps et avec les
applications r\'egulateurs $\R P_n, \R H_n$. Les deux complexes
sont en fait conjectur\'es \^etre le complexe motivique
$\underline{\Q(n)}$ sur $ {\mathrm Spec} E$.

Si on demanderait une application $r_{{\cal H}_n {\cal P}_n}:
{\cal H}_n(E) \to {\cal P}_n(E)$, la conjecture serait vraie pour
$n\leq 3$ mais fausse pour $n\geq 4$ (voir le th\'eor\`eme $4.7.$
de [G5]).

Les r\'esultats des paragraphes $\ref{2.1.},\ref{2.2.}$ reposent
sur le formalisme de la K-th\'eorie alg\'ebrique et du
r\'egulateur de Beilinson et on s'attend que des conjectures comme
la conjecture $\ref{conj3}$ r\'esultent du m\^eme type de
formalisme. Par contre, on s'attend que pour prouver la conjecture
$\ref{conj4}$ on aura \`a ma\^itriser une combinatoire tr\`es
difficile. On verra quelques exemples dans la section suivante.

\section{Le cas $n=4$}
Le hyperlogarithme $H(a_0|a_1,\cdots,a_n|a_{n+1})$ est invariant
par les transformations $a_i \to \alpha a_i + \beta, \alpha \in
{\C}^{\times},\beta \in \C$, donc il d\'epend en fait de quatre
variables. Le polylogarithme d\'epend d'une seule variable, donc
pour prouver la conjecture de Zagier dans le cas $n=4$ il faut
passer de quatre variables \`a une seule.

\subsection{Passer de quatre \`a deux variables}
Pour quatre nombres $A,B,C,D$ d'un corps $E$, on note par $ABCD$
leur birapport $\frac{(A-C)(B-D)}{(A-D)(B-C)}$. Pour deux nombres
$x\neq 0$ et $y\neq 1$ d'un corps $E$, on note
$[x,y]_{3,1}:=[0|x,0,0,y|1]$. Le th\'eor\`eme principal de cet
article est:
\begin{theor}
\label{theor3} Pour tout corps $E$, on a l'\'egalit\'e suivante
dans ${\cal H}_4(E)$:
\begin{\eq}
\label{9000} [A|B,C,D,E|F]=f(A,B,C,D,E)-f(B,C,D,E,F),
\end{\eq}
o\`u
$$-20f(A,B,C,D,E)=g(A,B,C,D,E)-g(\infty,B,C,D,E)-g(A,\infty,C,D,E)$$
$$-g(A,B,\infty,D,E)-g(A,B,C,\infty,E)-g(A,B,C,D,\infty)$$
$$-10{\mathrm cycl}[\frac{B-C}{B-A}]_4-10{\mathrm cycl}[\frac{A-B}{A-D}]_4+10{\mathrm cycl}[\frac{B-A}{B-D}]_4-10{\mathrm cycl}[\frac{D-B}{D-A}]_4,$$
o\`u
$$g(A,B,C,D,E)={\mathrm
cycl}\{[ABCD,BCDE]_{3,1}-[EDCB,EDCA]_{3,1}$$
$$-3[ABDC,ABDE]_{3,1}+3[EDBC,EDBA]_{3,1} \}.$$
Pour toute fonction $h$ de cinq variables $A_1,A_2,A_3,A_4,A_5$,
on a not\'e
$${\mathrm cycl} h(A_1,A_2,A_3,A_4,A_5)=\sum_{i=1}^5
h(A_i,A_{i+1},\cdots,A_5,A_1,\cdots,A_{i-1}).$$
\end{theor}

La fonction $f(A,B,C,D,E)$ peut \^etre interpr\'et\'e comme
$[A|B,C,D,E|\infty]$, i.e. comme r\'egularisation de l'int\'egrale
divergente $H(A|B,C,D,E|\infty)$.

\par{\bf Preuve: }
Une fois la r\'elation $(\ref{9000})$ d\'evin\'ee, la preuve est
un long calcul d'alg\'ebre lin\'eaire. En effet, on consid\`ere
les variables $A,B,tC+(1-t)A,tD+(1-t)A,E,F$ dans $E(t)$. D'apr\`es
la d\'efinition de ${\cal H}_4(E)$, il suffit de prouver que
$(\ref{9000})$ est tautologique en $t=0$ (calcul facile) et que
$\delta_4((\ref{9000}))$ est une \'egalit\'e dans ${\cal
H}_3(E(t)) \otimes E(t)_{\Q}^{\times}$. On peut continuer et
r\'eduire l'\'egalit\'e $\delta_4((\ref{9000}))$ \`a l'\'egalit\'e
$(\delta_3 \otimes id) \delta_4((\ref{9000}))$ et ensuite \`a
l'\'egalit\'e $(\delta_2 \otimes id \otimes id)(\delta_3 \otimes
id) \delta_4((\ref{9000}))$ dans $\Lambda^2 \fxq \otimes \fxq
\otimes \fxq$ pour $F=E(t,u,v), t,u,v$ variables.

\subsection{Passer de deux \`a une variables}
On note $\delta_{2,2}:{\cal H}_4(E)\to \Lambda^2 {\cal H}_2(E)$ la
composante de $\delta_4$ dans $\Lambda^2 {\cal H}_2(E)$. Le fait
que ${\cal H}_n(E) = {\cal P}_n(E)$ pour $n\leq 3$ et une chasse
au diagramme triviale sur la diagramme
$$\ \ \ \ {\cal P}_4(E)\stackrel{\delta_4}{\longrightarrow}{\cal P}_3(E)\otimes \exq$$
$$\downarrow\ \ \ \ \ \ \ \ \ \ \ \ \ \ \ \ \ \downarrow$$
$$\ \ \ \ \ \ \ \ \ \ \ \ \ \ \ \ \ \ \ {\cal H}_4(E)\stackrel{\delta_4}{\longrightarrow}{\cal H}_3(E)\otimes \exq \oplus \Lambda^2 {\cal H}_2(E)$$
montrent que la conjecture $\ref{conj4}$ est impliqu\'ee par la
conjecture

\begin{conj}
\label{conj5} Tout \'el\'ement de ${\cal H}_4(E)$ annul\'e par
$\delta_{2,2}$ provient de ${\cal P}_4(E)$.
\end{conj}

Il existe un \'el\'ement int\'eressant de ${\cal H}_4(E)$ annul\'e
par $\delta_{2,2}$. Soient $x\neq y$ et $z$ trois \'el\'ements de
$E\setminus\{ 0,1 \}$. On calcule $\delta_{2,2} [x,z]_{3,1} =
[x]_2 \wedge [z]_2$ dans $\Lambda^2 {\cal H}_2(E)=\Lambda^2 {\cal
P}_2(E)$. Il est classique que l'image de l'\'el\'ement
$$A(x,y)=[x]-[y]-[x/y]+[(1-x)/(1-y)]-[(1-1/x)/(1-1/y)]$$
de $\Q [E\setminus \{ 0,1 \} ]$ est nulle dans ${\cal P}_2(E)$.
Donc $\delta_{2,2} B(x,y;z) =0$, o\`u
$$B(x,y;z)=[x,z]_{3,1}-[y,z]_{3,1}-[x/y,z]_{3,1}+[(1-x)/(1-y),z]_{3,1}-[(1-1/x)/(1-1/y),z]_{3,1}.$$
Il est naturel de conjecturer ([G5])

\begin{conj}
\label{conj6} $B(x,y;z)\in {\cal P}_4(E)$, i.e. $B(x,y;z)$ est une
combinaison lin\'eaire de $4-$logarithmes.
\end{conj}
Cette conjecture \`a l'apparence anodine a r\'esist\'e aux
attaques de plusieurs math\'e\-ma\-ti\-ciens. En plus:

\begin{theor}
\label{theor4} La conjecture $\ref{conj3}$ pour $n=4$ et la
conjecture $\ref{conj6}$ impliquent la conjecture de Zagier pour
$n=4$.
\end{theor}

\par{\bf Preuve: }
On a vu que la conjecture de Zagier est impliqu\'ee par les
conjectures $\ref{conj3}$ et $\ref{conj4}$ et que la conjecture
$\ref{conj4}$ est impliqu\'ee par la conjecture $\ref{conj5}$. Il
suffit donc de prouver que la conjecture $\ref{conj6}$ implique la
conjecture $\ref{conj5}$.

Soit $R_2(E)\subset \Q [E\setminus \{ 0,1 \} ]$ le sous-espace
vectoriel engendr\'e par les \'el\'ements $A(x,y)$. La preuve de
la proposition $1.22.$ de [G1] et le fait que
$K_3^{[2]}(E)_{\Q}=K_3^{[2]}(E(t))_{\Q}$ montrent que ${\cal
R}_2^{\cal P}(E)=R_2(E)$.

On consid\`ere un \'el\'ement $H$ de ${\cal H}_4(E)$ annul\'e par
$\delta_{2,2}$. Par le th\'eor\`eme $\ref{theor3}$, il peut \^etre
\'ecrit comme $\sum_i a_i [x_i,z_i]_{3,1}+\sum_j b_j [y_j]_4$,
avec $a_i,b_j \in \Q$ et $x_i,z_i,y_j \in E$. On a donc $\sum a_i
[x_i]_2 \wedge [z_i]_2 =0$ dans $\Lambda^2 {\cal P}_2(E)$. Mais
${\cal P}_2(E)=\Q [E\setminus \{ 0,1 \} ]/R_2(E),$ donc $\sum a_i
[x_i] \otimes [z_i]$ est une combinaison lin\'eaire d'\'el\'ements
de type $[t]\otimes [u] + [u]\otimes [t]$ et de type $A(x,y)
\otimes [z]$ dans $\Q [E\setminus \{ 0,1 \} ] \otimes \Q
[E\setminus \{ 0,1 \} ]$. Donc $H$ est une combinaison lin\'eaire
d'\'el\'ements de type $[t,u]_{3,1}+[u,t]_{3,1}$ ($=0$, calcul
facile), $[y]_4$ et $B(x,y;z)$.

\vs \par{\bf Bibliographie:}

[BD]: A. A. Beilinson, P. Deligne: Interpr\'etation motivique de
la conjecture de Zagier reliant polylogarithmes et r\'egulateurs,
Proc. Sympos. Pure Math., vol. 55, Part 2, AMS, Providence, RI
(1994), p. 97-121

[B]: A. Borel: Cohomologie de $SL_n$ et valeurs de fonctions
z\^eta aux points entiers, Ann. Scuola Norm. Sup. Pisa Cl. Sci.
4(1977), p. 613-636

[G1]: A. B. Goncharov: Geometry of configurations, polylogarithmes
and motivic cohomology, Adv. Math. 114(1995), p. 197-318

[G2]: A. B. Goncharov: Explicit construction of characteristic
classes, Adv. Sov. Math. 16(1993), p. 169-210

[G3]: A. B. Goncharov: Chow polylogarithms and regulators, Math.
Res. Letters 2(1995), p. 99-114

[G4]: A. B. Goncharov: Multiple $\zeta$-numbers, hyperlogarithms
and mixed Tate motives, Preprint MSRI 058-93(1993)

[G5]: A. B. Goncharov: Polylogarithmes and motivic Galois group,
Proc. Sympos. Pure Math., vol. 55, Part 2, AMS, Providence, RI
(1994), p. 43-96

[Z1]: D. Zagier: Polylogarithms, Dedekind zeta functions and the
algebraic K-theory of fields, Progr. Math, vol. 89(1991), p.
391-430

[Z2]: D. Zagier: Hyperbolic manifolds and special values of
Dedekind zeta functions, Invent. Math. 83(1986), p. 285-301

\vs Institute of Mathematics of the Romanian Academy, Calea
Grivitei 21, 010702 Bucharest, Romania

E-mail: ndan@dnt.ro

\end{document}